\documentclass[11pt]{articleplus}   

\title{Singularities of Hinge Structures}
\author{Ciprian Borcea and Ileana Streinu}
\date{}

\begin{document}

\maketitle

\noindent
\begin{abstract}
Motivated by the hinge structure present in protein chains and other molecular conformations, we study the singularities of certain maps associated to body-and-hinge and panel-and-hinge chains. These are sequentially articulated systems where two consecutive rigid pieces are connected by a hinge, that is, a codimension two axis.

The singularities, or critical points, correspond to a dimensional drop in the linear span of the axes, regarded as points on a Grassmann variety in its Pl\"{u}cker embedding. These results are valid in arbitrary dimension. The three dimensional case is also relevant in robotics.

\end{abstract}

\section*{Introduction}
\label{sec:introduction}

\medskip \noindent
A {\em hinge} in the Euclidean space $R^d$ is formed when two $d$-dimensional bodies or two $(d-1)$-dimensional
panels are articulated along a common $(d-2)$-dimensional affine space (the hinge axis), so that the {\em
possible relative motions} of one object with respect to the other consist only of {\em rotations fixing the
given hinge axis}. Motion along the hinge axis is prohibited.

\begin{figure}[h]
	\vspace{-0.25 in}
    \centering
    {\includegraphics[width=12cm]{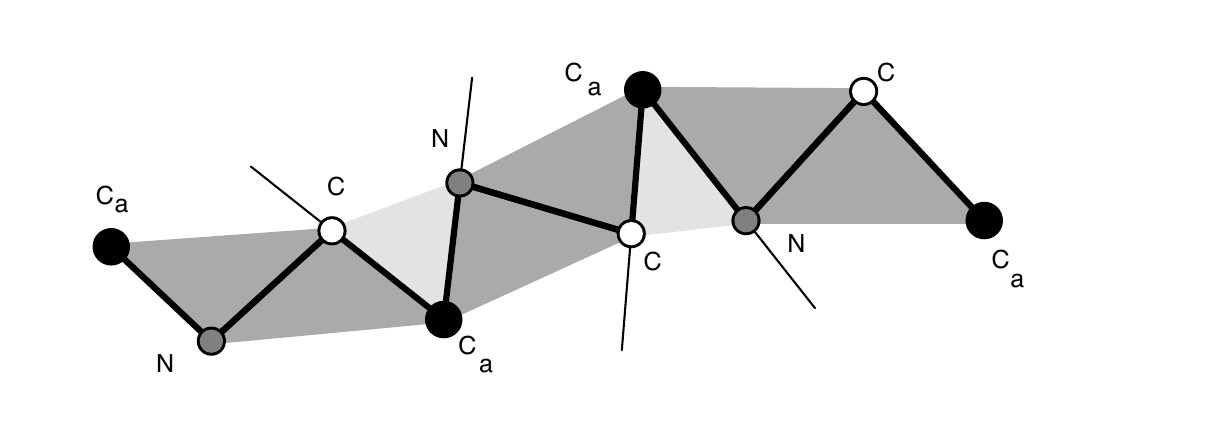}}\vspace{-0.3 in}
    \caption{A panel-and-hinge model for a protein backbone
	sequence. $C_{\alpha}$ atoms are represented by black dots, $N$
	atoms by grey and  $C$ atoms by white dots. The peptide planes
	containing $C_{\alpha} - C - N - C_{\alpha}$ bonds (dark grey)
	alternate with planes containing bonds $N - C_{\alpha} - C$ (light
	grey). The axes of the structure, shown as extended line segments,
	run along the $N - C_{\alpha}$ and $C_{\alpha} - C$ bond vectors.}
    \label{fig:peptide}
\end{figure}

\medskip \noindent
This situation appears for {\em molecular conformations} in $R^3$, when part of a molecule rotates with respect
to the remaining part around an {\em axis} corresponding to a chemical bond. Figure 1 schematically represents a
piece of a protein backbone \cite{BT} as a panel-and-hinge structure.

\medskip \noindent
We consider {\em ordered chains} of $n$ {\it bodies} or codimension-one {\it panels} $B_i, \ i=1, \cdots ,n$,
which are articulated {\em serially} by $n-1$ hinges $A_j, \ j=1,\cdots ,n-1$, with hinge $A_j$ linking $B_j$
and $B_{j+1}$. From now on, hinge axes will be simply called {\it hinges} or {\it axes} and will refer to the
corresponding codimension-two affine subspace of the chain configuration.

\begin{figure}[h]
	\vspace{-0.25 in}
    \centering
    {\includegraphics[width=14cm]{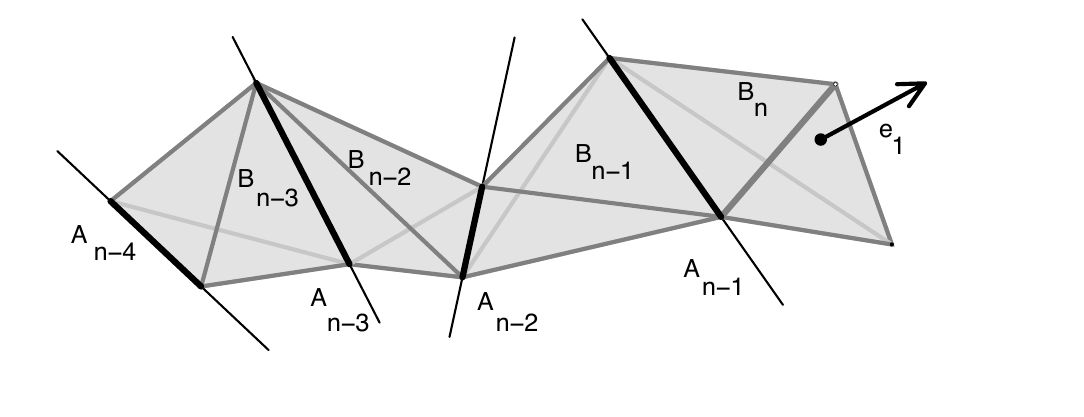}}\vspace{-0.3 in}
    \caption{End part of a body-and-hinge chain in $R^3$,
	with successive bodies identified as tetrahedra $B_j, B_{j+1}$,
	hinged along a common edge supported by the axis line $A_j$. The
	last body $B_n$ has a (rigidly) attached one-frame indicated by
	the vector $e_1$ (with given origin and direction relative to
	$B_n$).}
    \label{fig:bodyFrame}
\end{figure}

\medskip \noindent
We assume that our abstract objects (bodies or panels) can move through one another. By identifying
configurations which differ only by some rigid motion of the whole chain, the total {\em configuration space} is
naturally parametrized by the $(n-1)$-torus $(S^1)^{n-1}$. We factor out these rigid motions by {\em fixing the
first object}. This also fixes the first hinge. Clearly, a chain of hinged panels is simply a chain of hinged
bodies subject to the condition that two consecutive hinge axes span only a codimension-one affine subspace (the
corresponding panel).

\medskip \noindent
To the last object, we may attach some {\em frame} (e.g. a point, or a Cartesian $k$-frame) or some {\em flag}
(i.e. a sequence of linear subspaces, one included in the next), and study the {\em end-frame} or {\em end-flag
map} which takes a configuration to its corresponding frame or flag position. Note that the target is itself a
manifold (of frames or flags) and the resulting map is differentiable.

\medskip \noindent
We study the {\em singularities}, or {\em critical points} of such maps, that is configurations corresponding to a {\em drop in the rank of the differential}. We obtain geometrical characterizations of these singularities (Theorems \ref{thm:rank}, \ref{thm:plucker}, \ref{thm:confSpace}, \ref{thm:endFrame}) valid in arbitrary dimension $d$: they relate singular configurations to a lower dimensional span of the hinges in the corresponding projective Grassmann variety $G(d-1,d+1)\subset P_{{d+1\choose 2}-1}$. The most intuitive case, which was known in Robotics \cite{SDH,Bur1,Bur2}, is the end-point map in dimension 3, illustrated below.

\begin{theorem}
	\label{thm:endpoint3d}
Consider a body-and-hinge chain in $R^3$, with the first body fixed (i.e. identified with the ambient $R^3$) and with a marked point $e$ on the last body $B_n,\ n>3$. Consider the end-point map:

$$ e: (S^1)^{n-1} \rightarrow R^3\ , \ \ \theta =(\theta_1,...,\theta_{n-1})\mapsto e(\theta) $$

\noindent which registers, for a given configuration $\theta$ of the chain, the corresponding position
$e(\theta)$ of the marked point in the ambient space $R^3$.

\medskip \noindent
Then, the the differential of this map:\  $de(\theta)$ is of rank $<3$ if and only if there's a line through the end-point $e(\theta)$ which is projectively incident with all the axes $A_i(\theta),\ i=1,...,n-1$ of the corresponding configuration.
\end{theorem}

\medskip \noindent
{\em Projectively incident} means intersecting in $R^3$ or parallel, that is: intersecting `at infinity' in the projective completion $P_3=R^3\sqcup P_2$. 

\medskip \noindent
It should be emphasized that the intervention of a {\em projective characterization of singularities} is no accident - indeed, it echoes the known ``projective invariance of infinitesimal rigidity" in kinematics. See e.g. \cite{Wun} \cite{Weg}.

\medskip \noindent
We remark that, in dimension two, a body-and-hinge chain is as much as a panel-and-hinge chain, namely: a linkage given by $n$ rigid bars connected serially by revolute joints. This is, in other words, a {\em planar robot arm} and the singularities of the end-point map are known to be precisely the configurations with all bars along the same line \cite{Ha} \cite{KM1}, which indeed is the content of our result in dimension two. Thus, our hinge structures may be envisaged as higher dimensional versions of simple planar linkages. Although there is a conversion dictionary between a hinge-structure description and a linkage description in arbitrary dimension - as we outline in Section \ref{sec:cyclesToLinkages}, the former language seems better suited for characterizing singularities. We reinforce this aspect by discussing in Section \ref{sec:platforms} a related case in kinematics: infinitesimally flexible platforms.

\medskip
\noindent
The results in this paper have been presented at the Eighth International Symposium on {\em Effective Methods in Algebraic Geometry (MEGA)} 2005, Porto Conte, Alghero, Sardinia, May 26-June 2, 2005.

\section{The end-point map for body-and-hinge chains in $R^d$}
\label{sec:endpoint}

\noindent Let $B_1,...,B_n$ denote $d$-dimensional bodies in $R^d$. To be precise, one should think of each
$B_i$ as a copy of $R^d$, free to move relative to the ambient $R^d$. One may attach a Cartesian frame to the
copy and represent the movement of the body as the movement of the frame.

\medskip \noindent
We put a hinge $A_j$ between $B_j$ and $B_{j+1},\ j=1,...,n-1$, that is: we distinguish a codimension-two affine
subspace ( an axis) in $B_j$ and one in $B_{j+1}$, and we specify an isometry between them, and the two linked
bodies are now supposed to be positioned in the ambient $R^d$ subject to the condition that the two marked axes
coincide, and realize the specified isometry. The common axis, as seen in the ambient $R^d$, or on each of the
bodies so linked, will be denoted $A_j$.

\medskip \noindent
{\bf We shall identify the first body with the ambient $R^d$}, i.e. fix it as the reference body, because we are
interested in configurations only up to a rigid motion of the assembled chain.

\medskip \noindent
Clearly a hinge between two bodies allows one to rotate with respect to
the other, with the hinge axis remaining pointwise fixed. This relative
motion is parametrized by the unit circle $S^1$. Thus, with $B_1$ fixed,
the {\em configuration space} of the chain of $n$ hinged bodies is
parametrized by $(S^1)^{n-1}$.

\medskip \noindent
We distinguish now some particular point of the {\em last} body in the chain:\  $e\in B_n-A_{n-1}$, and call it
{\em the end-point}. (Obviously, we may assume $e$ to be away from the last axis $A_{n-1}$, since otherwise we
would restrict considerations to $n-1$ bodies.) This produces a map (to be denoted by $e$ as well):

$$ e: (S^1)^{n-1} \rightarrow R^d \ , \ \
\theta=(\theta_1,...,\theta_{n-1}) \mapsto e(\theta) $$

\noindent
which associates to a configuration $\theta$, the position $e(\theta)$
of the end-point with respect to the ambient space i.e. $B_1$. This
will be our {\em end-point map}.

\medskip \noindent
Our first concern is to describe the {\em singularities of the end-point
map}, that is: the configurations $\theta\in (S^1)^{n-1}$ where the
tangent map $de(\theta)$ has rank {\em strictly less} than its generic
rank. We have:

\begin{theorem}
	\label{thm:rank}
$rk(de(\theta))< d$ if and only if there's a line through $e(\theta)$ which is projectively incident with all
the axes $A_i(\theta), i=1,...,n-1$ of the corresponding configuration.
\end{theorem}

\medskip \noindent
Note that $A_1(\theta)=A_1$ is fixed, and the line through the end-point
in the theorem is either intersecting an axis or parallel to it
(i.e. meeting it ``at infinity'', when we complete $R^d$ to the
projective space $P_d$).

\medskip \noindent
{\em Proof:}\ The image of the differential $de(\theta)$ is spanned
by the tangent vectors at $e(\theta)$ to the circles
(or circles degenerated to a point) described by the end-point
$e(\theta)$ in the ambient $R^d$, when rotated around each axis
$A_i(\theta)$.

\medskip \noindent
This span is less than the full tangent space $R^d$ at $e(\theta)$ if and only if there's a line $\nu$ through
$e(\theta)$, normal to it. But $\nu$, will then be projectively incident with all axes.

\medskip \noindent
Indeed, if $e(\theta)$ happens to be on some axis, there's nothing to prove for that axis, while otherwise,
$\nu$ must lie in the hyperplane spanned by the axis under consideration, say $A_k(\theta)$ and $e(\theta)$,
which is the hyperplane normal to to the tangent at $e(\theta)$  for the circle described while rotating around
$A_k(\theta)$. This is, essentially, a partial derivative at $\theta$.

\medskip \noindent
By the same elementary theorem, if a line $\nu$ passing through $e(\theta)$
is projectively incident with all axes, it will be normal to
$im(de(\theta))$. \ \ \ $\Box$

\begin{corollary}
	\label{cor:orthogonal}
The space orthogonal to $im(de(\theta))$ is swept by all lines through
the end-point $e(\theta)$ which are projectively incident to all
axes $A_i(\theta), \ i=1,...,n-1$.
\end{corollary}

\begin{corollary}
	\label{cor:onto}
For $n> d$ and a generic choice of hinge axes, the differential of the
end-point map is generically onto, and its singularities are precisely
the configurations which allow some line through the end-point to be
projectively incident with all axes.
\end{corollary}

\medskip \noindent
{\bf Remarks:}\ i)\ The geometric argument used above does not even
require to be specific about the parametrization of the configuration
space by $(S^1)^{n-1}$, e.g. what position is considered for
$\theta=(0,..,0)$. It is enough to follow the infinitesimal displacements
of the end-point resulting from rotating as one body the part of the
chain from $B_i$ on, around $A_{i-1}$.

\medskip \noindent
ii)\ This approach also shows that {\em for infinitesimal considerations},
the {\em order of the axes} may turn out to be irrelevant, while clearly
essential otherwise.

\medskip \noindent
iii)\ We have emphasized in our statements the {\em purely projective characterization of the singularities}.
This is consistent with, in fact tantamount to the related phenomenon for linkages (cf. the so-called
``projective invariance of infinitesimal rigidity'' \cite{Weg}).

\medskip \noindent
For {\bf chains of hinged panels in $R^3$}, the line $\nu$ in the theorem must either pass through the
projective intersection of the two axes of an intermediate panel, or be contained in it; and is always contained
in the last panel.

\section{$k$-frames in $R^d$ and end-frame maps}
\label{sec:kframes}

\medskip \noindent
We begin by reviewing a few facts about the {\em homogeneous manifolds}
$W(k,d)$ defined by {\em all orthogonal (i.e. Cartesian) $k$-frames in
$R^d$}.

\medskip \noindent
One such frame consists of a point in $R^d$ (to be thought of as the origin of the frame) and $k$ ordered unit
vectors which are mutually orthogonal.

\medskip \noindent
Clearly, for $k=d$, we can identify the manifold $W(d,d)$ of all $d$-frames in
$R^d$ with the {\em isometry group} $Isom(R^d)$ of $R^d$:

$$ W(d,d)= R^d\rtimes O_R(d)\approx Isom(R^d) $$

\noindent
where $O_R(d)$ stands for the real orthogonal group in dimension $d$,
consisting of all orthogonal $d\times d$ matrices. (The columns of an
orthogonal matrix are the vectors of a $d$-frame.) The pair $(t,M)$
gives the isometry: \ $x\mapsto Mx+t$.

\medskip \noindent
Suppose now $0\leq k\leq d$, and note that one can parametrize all systems of $k$ ordered, mutually orthogonal
unit vectors in $R^d$ by the homogeneous space $O_R(d)/O_R(d-k)$ (where $O_R(d-k)$ is identified with the
subgroup of $O_R(d)$ fixing the first $k$ vectors in the standard basis). (These homogeneous spaces are called
Stiefel manifolds.) This gives the general description:

$$ W(k,d)=R^d\rtimes O_R(d)/O_R(d-k)\approx Isom(R^d)/O_R(d-k) $$

$$  dim_R W(k,d)= d+{d\choose 2}-{{d-k}\choose 2}=
{{d+1}\choose 2}-{{d-k}\choose 2} $$

\medskip \noindent
Notice that there's a {\em natural action} of the group of isometries in dimension $k$: \ $Isom(R^k)$, on the
space of $k$-frames in $R^d$: \ $W(k,d)$.

$$ Isom(R^k)\times W(k,d) \rightarrow W(k,d) $$

\medskip \noindent A $k$-frame gives an identification of its span with $R^k$, and $Isom(R^k)$ acts by displacing
the given frame to the image of the standard basis. Thus the action preserves the $k$-plane spanned by a frame
(through its origin), that is: a $k$-frame and its transforms have the same `supporting' $k$-plane.

\medskip \noindent
Let us fix a $k$-frame $E_k$ in the last body $B_n$ of a hinged chain.
As in the case of a point $e=E_0$, this gives an {\bf end-frame map}:

$$  (S^1)^{n-1}\rightarrow W(k,d), \ \ \theta \mapsto E_k(\theta) $$

\noindent
takig a configuration $\theta$ to the corresponding position of
the end $k$-frame in the ambient $R^d$.

\medskip \noindent
Again, for the singularities of the end-frame map, only the positions
of the axes matter, not their ordering, and the remarks on the action
of $Isom(R^k)$ on $W(k,d)$ give:

\begin{proposition}
	\label{prop:singularities}
The singularities of the end-frame map:

$$ E_k: (S^1)^{n-1}\rightarrow W(k,d), \ \ \theta \mapsto E_k(\theta) $$

\noindent
depend only on the set of axes and the $k$-plane spanned by the
end-frame.
\end{proposition}

\medskip \noindent
It may be useful in this context to cosider explicitly the map which
takes a $k$-frame to the (affine) $k$-plane it spans in $R^d$, as a map
to the Grassmann variety $G(k+1,d+1)$ parametrizing all $(k+1)$ linear
subspaces in $R^{d+1}$, that is: all projective $k$-planes in $P_d$:

$$ \pi_k: W(k,d) \rightarrow G(k+1,d+1) $$

\medskip \noindent
The axes themselves can be seen as points in $G(d-1,d+1)\approx G(2,d+1)$,
and our proposition says that the singularities of the end-frame map
depend only on $A_i(\theta), i=1,...,n-1$, and $\pi_k(E_k(\theta))$
as a point of $(G(d-1,d+1)^{n-1}/{\cal S}_{n-1})\times G(k+1,d+1)$.

\medskip \noindent
In order to see what kind of geometrical characterization of singularities
should emerge, we look in the next section at the case $k=d-2$.
A $(d-2)$-frame on the last body may be interpreted as a `loose hinge',
or half a hinge, and if we prescribe its matching half i.e. a
$(d-2)$-frame in the ambient $R^d$ (which is the first body), we may
interpret {\em the fibers of the end-frame map} $E_{d-2}$, as
{\em configuration spaces of cycles of $n$ hinged bodies}, for various
placements of the closing hinge.

\section{From chains to cycles}
\label{sec:chainTocycles}

\medskip \noindent
When we take $k=d-2$, Proposition 2.1 says that the singularities
of $E_{d-2}$ depend on corresponding configurations of $n$ points
in the Grassmann variety $G(d-1,d+1)$. Since the singularities of
the map indicate singularities of the fibers, and the fibers, in this case,
are cycles of $n$ bodies with $n$ hinges, we see that the {\em order}
of the $n$ points in the Grassmannian is not actually relevant.

\medskip \noindent
We should add the remark that our set-up generalizes the case $d=2$ of
the planar `robot arm' and planar polygon spaces \cite{Ha} \cite{KM1}
\cite{Bor2}.
In that case, one has singularities if and only if all $n$ axes
(which are simply points in $R^2\subset P_2$) are {\em collinear}.
In general, we have:

\begin{theorem}
	\label{thm:plucker}
Suppose $n\geq {{d+1}\choose 2}$. Consider the Pl\"{u}cker
embedding of the Grassmann variety:

$$ G(d-1,d+1)\hookrightarrow P_{{{d+1}\choose 2}-1} $$

\noindent
The end-frame map for a chain of $n$ hinged bodies in $R^d$:

$$ E_{d-2}: (S^1)^{n-1}\rightarrow W(d-2,d) $$

\noindent
has a singularity at $\theta\in (S^1)^{n-1}$ if and only if the
$n$ points of $G(d-1,d+1)$ corresponding to the axes $A_i(\theta)$,
$i=1,...,n-1$ and the span $\pi_{d-2}E_{d-2}(\theta)$ of the end-frame,
all lie in some hyperplane
section of the Grassmannian (in its Pl\"{u}cker embedding).
\end{theorem}

\medskip \noindent
Note that: \ $dim_R W(d-2,d)={{d+1}\choose 2}-1$.

\medskip \noindent
In terms of cycles, we have the simpler, but equivalent formulation:

\begin{theorem}
	\label{thm:confSpace}
The configuration space parametrizing the possible positions
(up to Euclidean motions) of a
cycle of $n\geq {{d+1}\choose 2}$ hinged bodies in $R^d$
is singular whenever the $n$ axes, as points in

$$ G(d-1,d+1)\hookrightarrow P_{{{d+1}\choose 2}-1} $$

\noindent
span less than the whole ambient projective space of the
Grassmannian.
\end{theorem}

\medskip \noindent
Note that a generic (initial) position of the $n$ axes gives a configuration
space of dimension $n-{{d+1}\choose 2}$.

\medskip \noindent
{\em Proof:}\ We extend the argument presented by Bricard in Tome II, Note H of \cite{Br2}.

\medskip \noindent
An infinitesimal motion of our chain of $n$ hinged bodies in $R^d$
corresponds with relative infinitesimal motions for each couple:

$$ (B_2/B_1)\ (B_3/B_2)\ ... (B_n/B_{n-1})(B_1/B_n) $$

\noindent
which are all tangent to uniform rotations with axes $A_1,...,A_n$.

\medskip \noindent
A simple way to encode a uniform rotation around a codimension two axis
$A_i\subset R^d$
uses an arbitrary point $M_i\in A_i$ and an element
$\omega_i=v^i_1\wedge v^i_2\wedge ...\wedge v^i_{d-2} \in  \bigwedge^{d-2}
R^d$, where $v^i_1,...,v^i_{d-2}$ is a basis of the subspace $A_i-M_i$,
whose exterior power $\omega_i$ represents the angular velocity of
the rotation.

\medskip \noindent
The information $(M_i,\omega_i)$, which generalizes the notion of
sliding vector (vecteur glissant) in dimension three, can also be
presented as an exterior vector:

$$ T_i=0_dM_i\wedge \omega_i + e_{d+1}\wedge \omega_i =
0_{d+1}M_i\wedge \omega_i \in \bigwedge^{d-1} R^{d+1} $$

\noindent
when we consider $R^d$ as the affine subspace $x_{d+1}=1$ in $R^{d+1}$
with the origin $0_d=(0,...,0,1)$, so that $e_{d+1}=0_{d+1}0_d$.
Thus uniform rotations become
representatives for points in $G(d-1,d+1)$ determined by their (affine)
axes.

\medskip \noindent
The component $0_dM_i\wedge \omega_i$ in $T_i$ expresses the velocity of
$0_d\in R^d$ rotating with respect to $A_i$.

\medskip \noindent
When we fix representatives $\alpha_i\in \bigwedge^{d-1}R^{d+1}$ for
all axes $A_i\in G(d-1,d+1)$, we have: \ $T_i=\tau_i\alpha_i$.

\medskip \noindent
The result of the relative infinitesimal motions given by $(M_i,\omega_i),\
i=1,...,n$ on the corresponding couples $(B_{i+1}/B_i)$ is clearly the
identity when considered relative to one and the same body, say $(B_1/B_1)$.
Thus, generalizing the null torsor condition in dimension three, we must
have:

$$  \sum_{i=1}^n \tau_i\alpha_i =0 \ \ \ \ \ \ \ \ \ (T) $$

\medskip \noindent
Indeed, the resulting velocities must be zero at the origin $0_d$ and
elswhere:

$$ \sum_{i=1}^n 0_dM_i\wedge \omega_i =0 $$

$$ \sum_{i=1}^n PM_i\wedge \omega_i =0 \ \ \ \  i.e. \ \ \ \
 \sum_{i=1}^n P0_d\wedge \omega_i=0 \ \mbox{for any} \ P\in R^d $$

\medskip \noindent
The last condition gives:

$$ \sum_{i=1}^n \omega_i =0 $$

\noindent
and $(T)$ follows.

\medskip \noindent
The dimension of the space of solutions $(\tau_i)_i$ of equation $(T)$
is $n-rank(\alpha_i)_i$, hence: the configuration space
has a singularity if and only if
the axes span less than the whole ambient space of the Grassmannian.
\ \ \ $\Box$

\medskip \noindent
{\bf Remarks:}\ In the generically {\em rigid} case for cycles, namely  $n={{d+1}\choose 2}$, a singularity in
the configuration space means {\em infinitesimal flexibility}.

\medskip \noindent
In space $(d=3)$, we would have a cycle of 6 hinged bodies or panels. The case of 6 panels corresponds to the
{\em cyclo-hexane molecule}, and when phrased in terms of linkages to 1-skeleta of octahedra. Thus, our result
recovers characterizations of infinitesimal flexibility for objects of some long-standing interest \cite{Br1},
\cite{Ben}, \cite{Br2}. A hyperplane section of the Grassmann-Pl\"{u}cker quadric $G(2,4)\subset P_5$ is also
called a linear complex. \noindent A note of Darboux to Koenigs' {\em `Le\c{c}ons de cin\'{e}matique'}, p.431,
mentions a fact known to Chasles: {\em a twisted cubic in $P_3$  has all its tangents in the same linear complex
i.e. a rational normal cubic has all its tangents in the same hyperplane section of the Grassmannian $G(2,4)$}.

\medskip \noindent
Other examples of six lines in a linear complex come from Bricard's flexible octahedra. In particular, as
observed in \cite{Ben} (sect. 17), a six-cycle in $R^3$ with {\em hinges symmetric in pairs relative to an
axis}, has one degree of freedom of motion. To see that the hinges are linearly dependent in $G(2,4)\subset
P_5$, note that the symmetry in a line, as a projective transformation $T$, can be given by a diagonal matrix
with two $+1$ and two $-1$ eigenvalues, hence inducing an involution with two $+1$ and four $-1$ eigenvalues on
$\bigwedge^2 R^4$. Thus $\ell_i+T\ell_i,\ i=1,2,3$ are dependent.

\section{End-frame and end-flag maps}
\label{sec:endFrame}

\medskip \noindent
Suppose $0\leq k < d-2$, and consider a $k$-frame attached to the
last body of a chain. Clearly, any extension of this $k$-frame to
a $k+r$-frame gives a factorization of the end-frame map $E_k$
through $E_{k+r}$:

$$ (S^1)^{n-1} \rightarrow W(k+r,d) \rightarrow W(k,d) $$

\noindent
The differential of the last arrow is surjective at all points,
and it follows that the singularities of $E_k$ are contained in
the singularities of $E_{k+r}$, for any extension of the end-frame.
This leads to:

\begin{theorem}
	\label{thm:endFrame}
Let $0\leq k \leq d-2$ and $n > dim_RW(k,d)={{d+1}\choose 2}-{{d-k}\choose 2}$.

\noindent
The end-frame map for a chain of hinged bodies in $R^d$:

$$ E_k: (S^1)^{n-1}\rightarrow W(k,d) $$

\noindent has a differential of rank {\em less than} $dim_RW(k,d)$ at $\theta\in (S^1)^{n-1}$ if and only if,
the $n-1$ points of the Grassmann variety $G(d-1,d+1)$ corresponding to the axes $A_i(\theta)$, $i=1,...,n-1$
and the locus made of $\pi_{d-2}E_{d-2}(\theta)$ for all extensions of the $k$-frame $E_k(\theta)$ to a
$(d-2)$-frame $E_{d-2}(\theta)$, are contained in some hyperplane section of the Grassmannian (in its
Pl\"{u}cker embedding).

\medskip \noindent
For a generic initial position $\theta=0$ of the axes and the
end-frame:

$$ dim_R E_k^{-1}(E_k(0))= n-{{d+1}\choose 2}+{{d-k}\choose 2}-1 $$

\end{theorem}

\medskip \noindent
This statement shows that one may replace frames with flags
(which is, in fact, the natural thing to do from the {\em complex}
point of view), and consider the singularities of the end-flag map:

$$ F_k : (S^1)^{n-1} \rightarrow W(k,d) \rightarrow Fl(k,d) $$

\noindent
obtained by composition with $ W(k,d) \rightarrow Fl(k,d) $,
which associates to an orthogonal frame $E_k=\{ e_0,e_1,...,e_k \}$ the
projective flag in $P_d(R)=R^d\cup P_{d-1}(R)$ made of subspaces
spanned by the first $m$ elements in the frame, with $0\leq m\leq k$.

\medskip \noindent
{\em Proof:} \ We use the shorter notation $Isom(R^k)=E(k)$. Lie algebras will be denoted with corresponding
small case letters. Using the homogeneous space description:

$$ W(k,d)=E(d)/O(d-k) $$

\noindent we may identify the tangent space to $W(k,d)$ at $E_k(\theta)$ with $e(d)/o(d-k)$.

\medskip \noindent As in our argument for Theorem 1., the image of the tangent map at $\theta$
is spanned by the $(n-1)$ tangent vectors corresponding to rotations around each axis (with the rest of the
chain imagined as rigid, from that axis on). These vectors are represented in $e(d)$ by the corresponding
infinitesimal rotations.

\medskip \noindent They do not span the whole tangent space $e(d)/o(d-k)$ precisely when there's a linear
functional on $e(d)$, vanishing on $o(d-k)$ and all the $(n-1)$ infinitesimal rotations.

\medskip \noindent The theorem then is simply the reading of this statement when converted via the {\em linear
isomorphism}\footnote{Intuitively, the linear isomorphism comes from $R^d$ regarded as a sphere $S^d$ of {\em
`infinite radius'}. Formally, this can be treated as a `contraction' in Lie group theory \cite{Se}.}

$$ e(d) \approx o(d+1) $$

\noindent and the natural identification of skew-symmetric two-forms (in $(d+1)$ variables) with $(d-1)$
exterior vectors: \ $o(d+1)\approx \wedge^{d-1}(R^{d+1})$. Indeed, an infinitesimal rotation around a
codimension-two axis corresponds precisely with the exterior vector representing the axis as a point of the
Grassmanniann $G(d-1,d+1)$.

\section{Converting cycles into linkages}
\label{sec:cyclesToLinkages}

\medskip \noindent
In this section we describe (canonical) procedures for associating
linkages with $2n$ vertices and $(2d-1)n$ edges in $R^d$ to {\em generic}
cycles of $n$ hinged bodies in $R^d, \ d\geq 3$. This association will
permit the identification of the cycle configuration space with
corresponding components of the linkage configuration space.

\medskip \noindent
The indices for axes and bodies should be understood cyclically i.e.
modulo $n$.

\medskip \noindent
We need to distinguish between the case of odd and even dimension.

\medskip \noindent
Suppose $d$ is {\bf odd}, that is: $d=2k+1$. In the generic case, all
intersections of $k$ consecutive axes are lines.

$$ l_i=A_i\cap A_{i+1}\cap ... \cap A_{i+k-1} $$

\noindent
One should regard $l_i$ as part of $A_i$ (and moving with it as the
cycle deforms into other configurations).

\medskip \noindent
We {\em choose} two points in general position
on each of these $n$ lines. This gives exactly $2k$ points on
each axis, and exactly $2k+2$ points on any pair of consecutive
axes, which corresponds to a body. Thus, the $d$-simplex generated
by the $2k+2$ points marks the body, and we take as edges in our
linkage all edges belonging to one of these $n$ simplices. A final
count gives $(2d-1)n$ edges.

\medskip \noindent
{\bf Remark:}\ There is, in fact, a {\em canonical} way to choose
two points on each of the above lines. Indeed,
for every pair of consecutive lines, there's a unique common
perpendicular incident to both, and this gives one point on each
line in the pair. In the end, one has two points on each line.

\medskip \noindent
Suppose now $d$ {\bf even}, that is: $d=2k$. In the generic case,
all intersections of $k$ consecutive axes are points:

$$ p_i=A_i\cap A_{i+1}\cap ... \cap A_{i+k-1} $$

\noindent
We consider these $n$ points $p_i$, together with $n$ points
{\em chosen generically}, one in each intersection of $k-1$
consecutive axes:

$$ q_i\in A_i\cap A_{i+1}\cap ... \cap A_{i+k-2} $$

\noindent
This gives exactly $2k-1$ points in each axis, and $2k+1$ points in
any pair of two consecutive axes. As in the odd case, this leads to
a linkage with $2n$ vertices and $(2d-1)n$ edges which is the
1-skeleton of a complex consisting of $n$ simplices of dimension
$d$ which share cyclically, one with the next, a $(d-2)$-face.

\medskip \noindent
{\bf Remark:}\ Again, the generic case allows for a {\em canonical
choice} of the points $q_i$. Indeed, one may define $q_i$ as the
orthogonal projection of $p_{i+1}$ on the plane
$A_i\cap A_{i+1}\cap ... \cap A_{i+k-2}$.

\medskip \noindent
For more definiteness, we recall the notions of configuration space
envisaged here for cycles of hinged bodies, respectively linkages.

\medskip \noindent
For {\bf cycles}, we consider an initial position of axes
$A_i=A_i(0)\subset R^d, \ i=1,...,n$. Every pair of consecutive
axes $(A_i,A_{i+1})$ belongs to a rigid body $B_{i+1}$,
understood as a copy of $R^d$. $B_{i+1}$ can
move relative to $B_i$ by rotating with respect to the common axis
$A_i$. We consider $B_1$ identified with the ambient $R^d$, fix a
$(d-2)$-frame  in $A_n$, and define the {\em configuration
space} ${\cal C}(A_1,...,A_n)$ for our initial position $\theta=0$ as
{\em the fiber of the end-frame map} $E_{d-2}$ over the
initial position $E_{d-2}(0)$. {\em (See sections 2 and 3 above.)}
In formulae:

$$ E_{d-2}: (S^1)^{n-1} \rightarrow W(d-2,d) $$

\noindent
with $E_{d-2}(0)$\   a \ $(d-2)$-frame in $A_n=A_n(0)$, and:

$$ {\cal C}(A_1,...,A_n)=E_{d-2}^{-1}(E_{d-2}(0)) $$

\medskip \noindent
The configuration space so defined is (up to canonical identifications)
independent of the end-frame chosen in $A_n$, or a cyclic permutation
of the indices.

\medskip \noindent
Our results in section 3 imply that, for a {\em generic} initial
position, the configuration space will be a smooth submanifold of
$(S^1)^{n-1}$, of dimension $n-{{d+1}\choose 2}$. However, this
submanifold might have several connected components.

\medskip \noindent
We turn now to {\bf linkages}. A linkage $L$ is a {\em weighted graph},
with weights indicating the length ascribed to each edge.
A {\em realization} of $L$ in $R^d$ is a map from the vertex set of $L$
to $R^d$, such that any two vertices defining an edge are placed at
the distance required by its weight. The {\em configuration space}
${\cal C}(L)={\cal C}(L,d)$ of $L$ in $R^d$ is the space of all
realizations of $L$ in $R^d$, modulo Euclidean motions
(i.e. orintation preserving isometries of $R^d$).
Cf.\cite{Bor1} \cite{B-S}.

\medskip \noindent
In order to emphasize the fact that our {\em canonical procedure}
for converting (generic) cycles into linkages is
{\em independent} of the representative chosen in describing the
configuration space ${\cal C}(A_1,...,A_n)$, one can use the
following:

\begin{proposition}
	\label{prop:twoConfig}
Let $A=(A_i)_i, A'=(A'_i)_i$ denote two configurations in
${\cal C}(A)={\cal C}(A_1,...,A_n)={\cal C}(A')$, with canonically
associated linkages $L(A)$ and $L(A')$. Then:

$$ L(A)=L(A') $$

\noindent
This means that their graphs can be identified, and
the length ascribed to corresponding edges is the same.
\end{proposition}

\medskip \noindent
{\em Proof:}\ The graphs, in our case, are clearly the same:
1-skeleta of identically labeled simplicial complexes. Thus,
one has to verify only the edge-length matching.

\medskip \noindent
To see this, we consider the simplices in the canonical
realization of $L(A)$ as markers for the $n$ bodies $B_i$. Imagining the
cycle A unhinged at $A_n$ (but with the $(d-2)$ simplicial face marked
on both $B_1$ and $B_n$ - the equivalent of a $(d-2)$-frame), there is
a continuous deformation of the chain  (by some trajectory in
$(S^1)^{n-1}$ linking A to A') which ends-up as A' by matching again
the two marked $(d-2)$ simplicial faces in $B_1$ and $B_n$. But this
restores perforce all the incidences (and orthogonalities) defining
the canonical realization of $L(A')$ as incidences (and orthogonalities)
of the moved simplices of $L(A)$.   \ \ \ $\Box$

\medskip \noindent
{\bf Note:} \ The argument shows a little more: the corresponding
simplices are not only congruent, but {\em realized with the same
orientation}. In fact, the linkage configuration space ${\cal C}(L(A))$
does contain realizations with one or the other orientation for
some of the simplices in the underlying complex, and our
inclusion ${\cal C}(A)\subset {\cal C}(L(A))$ covers only those
components where all orientations are as given in the canonical
realization associated to A. Denoting this image by ${\cal C}(L(A)_c)$,
we obtain a diagram:

\medskip
$\begin{array}{cclcc}
{\cal C}(A) & \subset & (S^1)^{n-1} & \subset &
G(d-1,d+1)^{n-2} \\
\downarrow & \ & \ & \ & \uparrow \\
{\cal C}(L(A)_c) & \ & \subset & \ & (P_d)^{2n-(d+1)}
\end{array} $

\medskip \noindent
where, considering the first body $B_1$ as fixed, and the corresponding
simplex with $(d+1)$ vertices in $L(A)$ fixed as well, the last
column records the axes $A_2(\theta),...,A_{n-1}(\theta)$, respectively
the remaining $2n-(d+1)$ vertices placed in $R^d\subset P_d$, and their
relation through a (generically defined) {\em rational map}.

\begin{lemma}
	\label{lemma:genericCycle}
Let $A$ be a generic $n$-cycle in $R^d$ with axes $A_1,...,A_n$.
For $A_{n+1}\neq A_n$ sufficiently close to $A_n$, we have:

$$ {\cal C}(A_1,...,A_{n+1})=S^1\times {\cal C}(A_1,...,A_n) $$
\end{lemma}

\medskip \noindent
This is the analogon of ``a small cut at a vertex'' for polygon spaces.
The proof, as in that case, amounts to observing that the fibers of
the end-frame map:

$$ E_{d-2}: (S^1)^{n-1} \rightarrow W(d-2,d) $$

\noindent
over a small neighborhood of $E_{d-2}(0)$ can be identified with
${\cal C}(A_1,...,A_n)=E_{d-2}^{-1}(E_{d-2}(0)$. \ \ \ $\Box$

\section{Cycle invariants: moduli}
\label{sec:moduli}

\noindent
When we look at $n$-cycles as points in $G(d-1,d+1)^n$
{\em modulo the diagonal action of the group of Euclidean motions in}
$R^d$, we have a parameter space of dimension
$2(d-1)n - {{d+1}\choose 2}$
containing cycle configuration spaces of generic dimension
$n-{{d+1}\choose 2}$. Thus, a parameter space for the cycle
configuration spaces, that is: {\em a moduli space},
should have dimension $(2d-3)n$. In other words, we expect
$(2d-3)n$ continuous parameters (also called {\em invariants} or
{\em moduli}), to characterize a generic configuration
space, at least up to a finite number of possibilities.

\medskip \noindent
{\bf Example:} \ For the planar case $(d=2)$, with $n$-cycles
understood as $n$-gons with prescribed edges, the obvious
invariants are the {\em edge lengths} themselves. The
{\em admissible edge-length-vectors} make-up a polyhedreal cone
in $R^n$ (with section a second hypersimplex). The various topological
types for planar polygon configuration spaces are then described in
terms of a subdivision into chambers of this cone. \cite{Bor2} \cite{Ha}
\cite{KM1} \cite{N}

\medskip \noindent
According to the previous section, the canonical linkage $L(A)$ associated to a generic cycle
$A=(A_1,...,A_n)\in {\cal C}(A)$ may be envisaged as an {\em invariant} of the cycle configuration space in
dimension $d\geq 3$. This gives, upfront, $(2d-1)n$ edge lengths, but we have a number of {\em orthogonalities}
in all canonical realizations, which make $2n$ distances dependent on the remaining $(2d-3)n$. Thus, one obtains
$(2d-3)n$ invariants. However, the cone they determine in $R^{(2d-3)n}$ is more complicated than in the planar
case.

\medskip \noindent
Obviously {\em rescaling} does not change the structure of the
configuration spaces, and we may replace the cone with a transversal
section of dimension $(2d-3)n-1$, corresponding to ratios of invariants.

\section{Platforms}
\label{sec:platforms}

\medskip \noindent
In this section we present a complementary result in kinematics, meant to emphasize the fact that {\em line
geometry} (or dually, axis geometry), that is: the use of Grassmann varieties $G(2,d+1)\approx G(d-1,d+1)$,
provides a natural context for singularity issues.

\medskip \noindent
Our example may be envisaged as a generalization of a theorem of Desargues, form the `perspective' of
`platforms'. For background on infinitesimal rigidity we refer to \cite{B-S} and \cite{Weg}.

\medskip \noindent
A {\bf platform} in $R^d$ consists of two (rigid) bodies connected by ${{d+1}\choose 2}$ rigid bars with ends
$p_{ij}$ on one body, respectively $q_{ij}$ on the other. For $d=2$ we connect the vertices of two triangles
with three bars, and the resulting framework is infinitesimally flexible precisely when the triangles are in
perspective for the given pairing of vertices: in other words, they produce a Desargues configuration. This
generalizes to:

\begin{proposition}
	\label{prop:platform}
A platform in dimension $d$ is infinitesimally flexible if and only if the ${{d+1}\choose 2}$ lines defined by
the connecting bars $p_{ij}q_{ij}$ lie in a hyperplane section of the Grassmannian $G(2,d+1)\subset
P_{{{d+1}\choose 2}-1}$.
\end{proposition}

\medskip \noindent
{\em Proof;}\ We use a `projective' version of the paltform, by imagining $R^d$ as $x_0=1$ in $R^{d+1}$ and the
origin linked by bars to $p_{ij}$ and $q_{ij}$. We may consider the first body fixed, and the infinitesimal
motion of the second given by an anti-symmetric $(d+1)\times (d+1)$ matrix $A$. An infinitesimal motion of the
platform requires:

$$ <Aq_{ij},p_{ij}-q_{ij}>=0 \ \ i.e. \ \ <Aq_{ij},p_{ij}>=0 $$

\noindent and this linear system in the unknowns $a_{ij}$ has a non-trivial solution precisely when the exterior
two-vectors: \ $p_{ij}\wedge q_{ij}$, $1\leq i<j\leq (d+1)$ are linearly dependent. \ \ \ $\Box$

\medskip \noindent
{\bf Remark:}\ For the `Stewart-Gough platform' i.e. $d=3$, this fact is presented in \cite{RV}, following
\cite{Mer}.

\section{Conclusions}
\label{sec:conclusion}

\medskip \noindent
In principle, understanding the singularities of a map offers a key towards the topology of its fibers (usually
via some Morse theory), and this approach is common to a number of studies on mechanical linkages \cite{N}
\cite{Ha} \cite{KM1} \cite{KM2} \cite{Bor2}.

\medskip \noindent
One should also remark that, with configuration spaces of this nature, all relevant maps are {\em algebraic},
and this allows the intervention of complex algebraic-geometry, as in \cite{Kly} \cite{Bor2}.

\end{document}